\documentclass[a4paper,12pt]{amsart}
\usepackage{amsfonts}
\usepackage{amssymb}
\usepackage{ifthen}
\usepackage{graphicx}
\nonstopmode \numberwithin{equation}{section}
\usepackage[usenames]{color}
\newtheorem{thm}{Theorem}[section]
\newtheorem{lem}{Lemma}[section]
\newtheorem{cor}[thm]{Corollary}
\newtheorem{prop}[thm]{Proposition}

\newtheorem{step}{Step}[section]

\theoremstyle{definition}
\newtheorem{mlem}{Main lemma}[section]
\newtheorem{assertion}{Assertion}[section]
\newtheorem{cl}{Claim}[section]
\newtheorem{ca}{Case}[section]
\newtheorem{sca}{Subcase}[section]
\newtheorem{scl}{Subclaim}[section]
\newtheorem{conj}[thm]{Conjecture}
\newtheorem{fact}{Fact}[section]
\newtheorem{defn}[thm]{Definition}
\newtheorem{op}[thm]{Open Problem}
\newtheorem{ques}[thm]{Question}
\newtheorem{rem}[thm]{Remark}
\newtheorem{exam}[thm]{Example}

\numberwithin{equation}{section}

\newcounter {own}
\def\theown {\thesection       .\arabic{own}}

\newenvironment{pf}[1][]{%
 \vskip 3mm
 \noindent
 \ifthenelse{\equal{#1}{}}%
  {{\slshape Proof. }}%
  {{\slshape #1.} }%
 }%
{\qed\bigskip}

\newcounter{alphabet}
\newcounter{tmp}
\newenvironment{Thm}[1][]{\refstepcounter{alphabet}%
\bigskip%
\noindent%
{\bf Theorem \Alph{alphabet}}%
\ifthenelse{\equal{#1}{}}{}{ (#1)}%
{\bf .} \itshape}{\vskip 8pt}

\makeatletter
\newcommand{\Ref}[1]{\@ifundefined{r@#1}{}{\setcounter{tmp}{\ref{#1}}\Alph{tmp}}}
\makeatother

\newenvironment{Lem}[1][]{\refstepcounter{alphabet}%
\bigskip%
\noindent%
{\bf Lemma \Alph{alphabet}}%
{\bf .} \itshape}{\vskip 8pt}

\newcounter{alphabet2}

\newcommand{\ID}{{\mathbb D}}
\newcommand{\IB}{{\mathbb B}}

\def\be{\begin{equation}}
\def\ee{\end{equation}}

\newcommand{\ben}{\begin{enumerate}}
\newcommand{\een}{\end{enumerate}}

\newcommand{\blem}{\begin{lem}}
\newcommand{\elem}{\end{lem}}
\newcommand{\bthm}{\begin{thm}}
\newcommand{\ethm}{\end{thm}}
\newcommand{\bcor}{\begin{cor}}
\newcommand{\ecor}{\end{cor}}
\newcommand{\beg}{\begin{exam}}
\newcommand{\eeg}{\end{exam}}
\newcommand{\begs}{\begin{examples}}
\newcommand{\eegs}{\end{examples}}
\newcommand{\bdefe}{\begin{defn}}
\newcommand{\edefe}{\end{defn}}
\newcommand{\bprob}{\begin{prob}}
\newcommand{\eprob}{\end{prob}}
\newcommand{\bques}{\begin{ques}}
\newcommand{\eques}{\end{ques}}
\newcommand{\bei}{\begin{itemize}}
\newcommand{\eei}{\end{itemize}}
\newcommand{\bcon}{\begin{conj}}
\newcommand{\econ}{\end{conj}}
\newcommand{\bop}{\begin{op}}
\newcommand{\eop}{\end{op}}

\newcommand{\bas}{\begin{assertion}}
\newcommand{\eas}{\end{assertion}}

\newcommand{\bfa}{\begin{fact}}
\newcommand{\efa}{\end{fact}}

\newcommand{\bca}{\begin{ca}}
\newcommand{\eca}{\end{ca}}

\newcommand{\bst}{\begin{step}}
\newcommand{\est}{\end{step}}

\newcommand{\bsca}{\begin{sca}}
\newcommand{\esca}{\end{sca}}

\newcommand{\bcl}{\begin{cl}}
\newcommand{\ecl}{\end{cl}}

\newcommand{\bmlem}{\begin{mlem}}
\newcommand{\emlem}{\end{mlem}}

\newcommand{\bscl}{\begin{scl}}
\newcommand{\escl}{\end{scl}}

\newcommand{\bcons}{\begin{conjs}}
\newcommand{\econs}{\end{conjs}}

\newcommand{\bprop}{\begin{prop}}
\newcommand{\eprop}{\end{prop}}

\newcommand{\br}{\begin{rem}}
\newcommand{\er}{\end{rem}}
\newcommand{\brs}{\begin{rems}}
\newcommand{\ers}{\end{rems}}
\newcommand{\bo}{\begin{obser}}
\newcommand{\eo}{\end{obser}}
\newcommand{\bos}{\begin{obsers}}
\newcommand{\eos}{\end{obsers}}
\newcommand{\bpf}{\begin{pf}}
\newcommand{\epf}{\end{pf}}
\newcommand{\ba}{\begin{array}}
\newcommand{\ea}{\end{array}}
\newcommand{\beq}{\begin{eqnarray}}
\newcommand{\beqq}{\begin{eqnarray*}}
\newcommand{\eeq}{\end{eqnarray}}
\newcommand{\eeqq}{\end{eqnarray*}}

\newcommand{\ds}{\displaystyle}

\newcounter{minutes}
\divide\time by 60
\newcounter{hours}
\multiply\time by 60 \addtocounter{minutes}{-\time}

\begin{document}

\bibliographystyle{amsplain}
\title []
{Schwarz type lemma, Landau type theorem and Lipschitz type
space of solutions to inhomogeneous biharmonic equations}

\def\thefootnote{}
\footnotetext{ \texttt{\tiny File:~\jobname .tex,
          printed: \number\day-\number\month-\number\year,
          \thehours.\ifnum\theminutes<10{0}\fi\theminutes}
} \makeatletter\def\thefootnote{\@arabic\c@footnote}\makeatother

\author{Shaolin Chen}
 \address{Sh. Chen, College of Mathematics and
Statistics, Hengyang Normal University, Hengyang, Hunan 421008,
People's Republic of China.} \email{mathechen@126.com}

\author{Peijin Li}
\address{P. Li, Department of Mathematics, Hunan First Normal University, Changsha,
Hunan 410205, People's Republic of China.} \email{wokeyi99@163.com}

\author{Xiantao Wang}
\address{X. Wang, Department of Mathematics, Shantou University, Shantou,
Guangdong 515063, People's Republic of China.}
\email{xtwang@stu.edu.cn}



\begin{abstract}
The purpose of this paper is to study the properties of the
solutions to the inhomogeneous biharmonic equations: $\Delta(\Delta f)=g$, where
$g:$ $\overline{\mathbb{D}}\rightarrow\mathbb{C}$ is a continuous
function and $\overline{\mathbb{D}}$ denotes the closure of the unit
disk $\mathbb{D}$ in the complex plane $\mathbb{C}$. In fact, we
establish the following properties for those solutions: Firstly, we
establish the Schwarz type lemma. Secondly, by using the obtained
results, we
 get a Landau type theorem. Thirdly, we discuss their Lipschitz type property.
\end{abstract}

\maketitle \pagestyle{myheadings} \markboth{ Shaolin Chen, Peijin Li and Xiantao
Wang}{Schwarz type lemma, Landau type theorem and Lipschitz type space}

\section{Preliminaries and  main results }\label{csw-sec1}

Let $\mathbb{C}$ denote the complex plane. For
$a\in\mathbb{C}$ and  $r>0$, we let $\ID(a,r)=\{z:\, |z-a|<r\}$, $\mathbb{D}_r=\mathbb{D}(0,r)$ and
$\mathbb{D}=\mathbb{D}_1$, the open unit disk in
$\mathbb{C}$. Let $\mathbb{T}=\partial\mathbb{D}$ be the boundary of
$\mathbb{D}$, and $\overline{\ID}=\ID\cup \mathbb{T}$, the closure of $\mathbb{D}$. Furthermore, we denote by $\mathcal{C}^{m}(\Omega)$ the set of all
complex-valued $m$-times continuously differentiable functions from
$\Omega$ into $\mathbb{C}$, where $\Omega$ stands for a subset of
$\mathbb{C}$ and $m\in\mathbb{N}_0:=\mathbb{N}\cup\{0\}$. In
particular, $\mathcal{C}(\Omega):=\mathcal{C}^{0}(\Omega)$ denotes the
set of all continuous functions in $\Omega$.

For a real $2\times2$ matrix $A$, we use the matrix norm
$$\|A\|=\sup\{|Az|:\,|z|=1\}$$ and the matrix function
$$\lambda(A)=\inf\{|Az|:\,|z|=1\}.$$

For $z=x+iy\in\mathbb{C}$, the
formal derivative of a complex-valued function $f=u+iv$ is given
by
$$D_{f}=\left(\begin{array}{cccc}
\ds u_{x}\;~~ u_{y}\\[2mm]
\ds v_{x}\;~~ v_{y}
\end{array}\right),
$$
so that
$$\|D_{f}\|=|f_{z}|+|f_{\overline{z}}| ~\mbox{ and }~ \lambda(D_{f})=\big| |f_{z}|-|f_{\overline{z}}|\big |,
$$
where $$f_{z}=\frac{1}{2}\big(
f_x-if_y\big)\;\;\mbox{and}\;\; f_{\overline{z}}=\frac{1}{2}\big(f_x+if_y\big).$$  We use
$$J_{f}:=\det D_{f} =|f_{z}|^{2}-|f_{\overline{z}}|^{2}
$$
to denote the {\it Jacobian} of $f$.

Let $\varphi, f^{\ast}\in \mathcal{C}(\mathbb{T})$,
$g\in\mathcal{C}(\overline{\mathbb{D}})$ and
$f\in\mathcal{C}^{4}(\mathbb{D})$. Of particular interest to us is
the following {\it inhomogeneous biharmonic equation} in
$\mathbb{D}$:

\be\label{eq-ch-1.0}\Delta(\Delta f)=g,\ee and its following
associated Dirichlet boundary value problem: \be\label{eq-ch-1}
\begin{cases}
\displaystyle  f_{\overline{z}}=\varphi &\mbox{ in }\, \mathbb{T},\\
\displaystyle f=f^{\ast}&\mbox{ in }\, \mathbb{T},
\end{cases}
\ee  where $$\Delta f=f_{x^2}+f_{y^2}=4f_{z \overline{z}}$$ is the
{\it Laplacian} of $f$. In particular, if $g\equiv0$, then the
solutions to (\ref{eq-ch-1.0}) are {\it biharmonic mappings} (see
\cite{AA,CWY,GMR-2018,S-2003}).

The inhomogeneous biharmonic equation arises in areas of continuum mechanics,
including linear elasticity theory and the solution of Stokes flows
(cf. \cite{Ha,Kh-1996, We}). Most important applications of the theory of functions of one complex variable were obtained in the plane theory of elasticity and in the
approximate theory of plates subject to normal loading (cf. \cite{HB-1965, L-1964}).
 Baernstein II and  Kovalev \cite{BK-2005} investigated
the H\"older continuity of the {\it gradient mapping}
$f\in\mathcal{C}^{1}(\Omega)$,
 where $\Omega$ is a subset of $\mathbb{C}$ and $f_{\overline{z}}=u$ for some $u\in\mathcal{C}(\Omega)$. In \cite{CP},
 the authors  studied the geometric properties
 of the gradient mappings.
 In this paper, we will discuss the behavior
and the potential properties of solutions to a class of
inhomogeneous biharmonic equations whose Dirichlet boundary values
are gradient mappings. This study continues the investigation in
\cite{BK-2005,CP} and is mainly motivated by the discussions in the
papers of  Abdulhadi and Abu Muhanna \cite{AA}, Colonna \cite{Co},
Heinz \cite{He}, Kalaj and Pavlovi\'c \cite{K2}, and the monograph
of Pavlovi\'c \cite{Pav1}. In order to state our main results, we
introduce some necessary terminologies.

For $z, w\in\mathbb{D}$, let
\be\label{G-1}G(z,w)=|z-w|^{2}\log\left|\frac{1-z\overline{w}}{z-w}\right|^{2}-(1-|z|^{2})(1-|w|^{2})\ee
and $$P(z,e^{i\theta})=\frac{1-|z|^{2}}{|1-ze^{-i\theta}|^{2}}$$
denote the {\it biharmonic Green function} and {\it (harmonic)
Poisson kernel}, respectively, where $\theta\in[0,2\pi]$.

By \cite[Theorem 2]{Be}, we see that all
solutions to the  inhomogeneous biharmonic equation (\ref{eq-ch-1.0}) satisfying the boundary value conditions
(\ref{eq-ch-1}) are given by

\beq\label{eq-ch-3}
f(z)&=&\mathcal{P}_{f^{\ast}}(z)+\frac{1}{2\pi}\int_{0}^{2\pi}\overline{z}e^{it}f^{\ast}(e^{it})\frac{1-|z|^{2}}{(1-\overline{z}e^{it})^{2}}dt\\
\nonumber
&&-(1-|z|^{2})\mathcal{P}_{\varphi_{1}}(z)-\frac{1}{16\pi}\int_{\mathbb{D}}g(w)G(z,w)dA(w),
\eeq where $dA(w)$ denotes the Lebesgue area measure in
$\mathbb{D}$,
$\varphi_{1}(e^{it})=\varphi(e^{it})e^{-it}$,
$$\mathcal{P}_{f^{\ast}}(z)=\frac{1}{2\pi}\int_{0}^{2\pi}P(z,e^{it})f(e^{it})dt$$ and
$$\mathcal{P}_{\varphi_{1}}(z)=\frac{1}{2\pi}\int_{0}^{2\pi}P(z,e^{it})\varphi_{1}(e^{it})dt.$$
In particular, if $f$ satisfies (\ref{eq-ch-3}) and is harmonic
(i.e. $\Delta f=0$) in $\mathbb{D}$, then for $z\in\mathbb{D}$,
\be\label{eq-clw-b-1}\mathcal{P}_{\varphi_{1}}(z)=\frac{1}{2\pi}\int_{0}^{2\pi}\frac{\overline{z}e^{it}f^{\ast}(e^{it})}{(1-\overline{z}e^{it})^{2}}dt~\mbox{and}~f(z)=\mathcal{P}_{f^{\ast}}(z).\ee
Moreover, if $f$ satisfies (\ref{eq-ch-3}) and is  biharmonic (i.e.
$\Delta(\Delta f)=0$) in $\mathbb{D}$, then for $z\in\mathbb{D}$,
$$f(z)=\mathcal{P}_{f^{\ast}}(z)+(1-|z|^{2})\left[\frac{1}{2\pi}\int_{0}^{2\pi}
\frac{\overline{z}e^{it}f^{\ast}(e^{it})}{(1-\overline{z}e^{it})^{2}}dt-\mathcal{P}_{\varphi_{1}}(z)\right].$$

In \cite{MM}, the solvability of the inhomogeneous biharmonic equations has been studied. In the following, we will investigate the Schwarz type lemmas, the Landau type theorem and the Lipschitz continuity of solutions
to the inhomogeneous biharmonic equation (\ref{eq-ch-1.0}) with the boundary value conditions
(\ref{eq-ch-1}).

The classical Schwarz lemma states that an analytic function
$f$ from $\mathbb{D}$ into itself with $f(0)=0$ satisfies
$|f(z)|\leq|z|$ for all $z\in\mathbb{D}$. It is well-known that the
Schwarz lemma has become a crucial theme in a lot of branches of
mathematical research for more than a hundred years to date.

Heinz \cite{He} proved the following result, which is called the
Schwarz lemma of  harmonic mappings: If $f$ is a harmonic mapping
from $\mathbb{D}$ into itself with $f(0)=0$, then for
$z\in\mathbb{D}$, \be\label{eqh}|f(z)|\leq\frac{4}{\pi}\arctan
|z|.\ee

Later, Pavlovi\'c \cite[Theorem 3.6.1]{Pav1} improved (\ref{eqh})
and obtained the following general form:
\be\label{eq-pav1}\left|f(z)-\frac{1-|z|^{2}}{1+|z|^{2}}f(0)\right|\leq\frac{4}{\pi}\arctan
|z|,~z\in\mathbb{D},\ee where $f$ is a  harmonic mapping from
$\mathbb{D}$ into itself. See \cite{CK, CP-0, KV, K3} for more
discussions in this line.

By analogy  with the inequality (\ref{eq-pav1}),  we obtain the following result.
\begin{thm}\label{thm-1}
Suppose that $g\in\mathcal{C}(\overline{\mathbb{D}})$, $\varphi\in
\mathcal{C}(\mathbb{T})$ and suppose that
$f\in\mathcal{C}^{4}(\mathbb{D})\cap\mathcal{C}(\overline{\mathbb{D}})$
 satisfies {\rm
(\ref{eq-ch-1.0})} and {\rm (\ref{eq-ch-1})}, where
$f^{\ast}=f|_{\mathbb{T}}$. Then for $z\in\overline{\mathbb{D}}$,
\beq\label{eq-ch-4}
\bigg|f(z)&-&\frac{1-|z|^{2}}{1+|z|^{2}}\mathcal{P}_{f^{\ast}}(0)+\frac{(1-|z|^{2})^{2}}{1+|z|^{2}}\mathcal{P}_{\varphi_{1}}(0)\bigg|\\
\nonumber &\leq&
\frac{4}{\pi}\|\mathcal{P}_{f^{\ast}}\|_{\infty}\arctan|z|
+\frac{4\|\mathcal{P}_{\varphi_{1}}\|_{\infty}}{\pi}(1-|z|^{2})\arctan|z|\\
\nonumber&&+
\|\mathcal{P}_{|f^{\ast}|}\|_{\infty}|z|+\frac{1}{64}\|g\|_{\infty}(1-|z|^{2})^{2},
\eeq where $$\|g\|_{\infty}=\sup_{z\in\mathbb{D}}\{|g(z)|\},~
\|\mathcal{P}_{f^{\ast}}\|_{\infty}=\sup_{z\in\mathbb{D}}\{|\mathcal{P}_{f^{\ast}}(z)|\},~
\|\mathcal{P}_{|f^{\ast}|}\|_{\infty}=\sup_{z\in\mathbb{D}}\{|\mathcal{P}_{|f^{\ast}|}(z)|\},$$
  $\|\mathcal{P}_{\varphi_{1}}\|_{\infty}=\sup_{z\in\mathbb{D}}\{|\mathcal{P}_{\varphi_{1}}(z)|\}$
  and $\varphi_{1}$ is defined in {\rm(\ref{eq-ch-3})}.
Moreover, if we take $g(z)\equiv M$, a positive constant, and
$$f(z)=\frac{M}{64}(1-|z|^{2})^{2},$$ then the inequality {\rm
(\ref{eq-ch-4})} is sharp.

\end{thm}

\begin{rem} \begin{enumerate}
\item Under the hypothesis of Theorem \ref{thm-1}, if $g\equiv0$, then $f$ is a biharmonic mapping, and
(\ref{eq-ch-4}) can be written in the following form
\beq\label{eq-clw-4.0}\nonumber
\bigg|f(z)-\frac{1-|z|^{2}}{1+|z|^{2}}\mathcal{P}_{f^{\ast}}(0)+\frac{(1-|z|^{2})^{2}}{1+|z|^{2}}\mathcal{P}_{\varphi_{1}}(0)\bigg|
&\leq& \frac{4}{\pi}\|\mathcal{P}_{f^{\ast}}\|_{\infty}\arctan|z|
+\|\mathcal{P}_{|f^{\ast}|}\|_{\infty}|z|\\
&&+\frac{4\|\mathcal{P}_{\varphi_{1}}\|_{\infty}}{\pi}(1-|z|^{2})\arctan|z|.
\eeq  Furthermore, the biharmonic mapping $f(z)=1-|z|^{2}$
$(z\in\overline{\mathbb{D}})$ shows that (\ref{eq-clw-4.0}) is sharp
in $\mathbb{T}$.

\item
Under the hypothesis of Theorem \ref{thm-1}, if $f$ is harmonic in
$\mathbb{D}$, and maps $\mathbb{D}$ into itself, then, by
(\ref{eq-clw-b-1}), we see that
$$\mathcal{P}_{\varphi_{1}}(z)-\frac{1}{2\pi}\int_{0}^{2\pi}\frac{\overline{z}e^{it}f^{\ast}(e^{it})}{(1-\overline{z}e^{it})^{2}}dt=0,$$
which implies that $\mathcal{P}_{\varphi_{1}}(0)=0$ and
$f(z)=\mathcal{P}_{f^{\ast}}(z)$. Hence, by(\ref{eq-pav1}),
 we
have
\beqq\bigg|f(z)-\frac{1-|z|^{2}}{1+|z|^{2}}\mathcal{P}_{f^{\ast}}(0)+\frac{(1-|z|^{2})^{2}}{1+|z|^{2}}\mathcal{P}_{\varphi_{1}}(0)\bigg|&=&
\left|f(z)-\frac{1-|z|^{2}}{1+|z|^{2}}\mathcal{P}_{f^{\ast}}(0)\right|\\
&\leq&\frac{4}{\pi}\arctan |z|\eeqq for $z\in\mathbb{D}$.
\end{enumerate}
\end{rem}

Let us recall  the classical boundary Schwarz lemma of analytic
functions, which is as follows.

\begin{Thm}{\rm (\cite{G})}\label{Thm-S} Suppose that $f$ is an analytic function from
$\mathbb{D}$ into itself. If $f(0)=0$ and $f$ is analytic at $z=1$
with $f(1)=1$, then $f'(1)\geq1$. Moreover, the inequality is sharp.
\end{Thm}

This useful result has attracted much attention and has been
generalized in various forms (see, e.g., \cite{BK, CK, CP-0, K3,
Kra, LT, WZh}). In the following, by applying Theorem \ref{thm-1},
we establish a Schwarz type lemma at the boundary for the solutions
to the inhomogeneous biharmonic equation \eqref{eq-ch-1.0}.

\begin{thm}\label{thm-2}
Suppose that $g\in\mathcal{C}(\overline{\mathbb{D}})$, $\varphi\in
\mathcal{C}(\mathbb{T})$, and suppose that
$f\in\mathcal{C}^{4}(\mathbb{D})$ satisfies the following equations:
$$\begin{cases}\displaystyle \Delta(\Delta f)=g &\mbox{ in }\, \mathbb{D},\\
\displaystyle  f_{\overline{z}}=\varphi &\mbox{ in }\, \mathbb{T},\\
\displaystyle f=\phi&\mbox{ in }\, \mathbb{T},
\end{cases}$$ where $\varphi_{1}$ is defined in {\rm
(\ref{eq-ch-3})},
  $\phi\in\mathcal{C}(\overline{\mathbb{D}})$ is analytic in $\mathbb{D}$, and, further, $\varphi_{1}$ and $\phi$ satisfy that
  $|\phi|\leq1,$
  $\|\mathcal{P}_{\varphi_{1}}\|_{\infty}<\frac{1}{2}$ and $$|\mathcal{P}_{\phi}(0)|<\frac{1-2\|\mathcal{P}_{\varphi_{1}}\|_{\infty}}{1+2\|\mathcal{P}_{\varphi_{1}}\|_{\infty}}.$$
If $\lim_{r\rightarrow1^{-}}|f(r\eta)|=1$ for $\eta\in\mathbb{T}$, then
\be\label{Sch-1}\liminf_{r\rightarrow1^{-}}\frac{|f(\eta)-f(r\eta)|}{1-r}\geq\frac{1-|\mathcal{P}_{\phi}(0)|}{1
+|\mathcal{P}_{\phi}(0)|}-2\|\mathcal{P}_{\varphi_{1}}\|_{\infty}.\ee
In particular, if $\|\mathcal{P}_{\varphi_{1}}\|_{\infty}=0$, then
the inequality of {\rm(\ref{Sch-1})} is sharp.
\end{thm}

Colonna \cite{Co} obtained a sharp Schwarz-Pick type lemma for
harmonic mappings, which is as follows: If $f$ is a harmonic mapping
from $\mathbb{D}$ into itself, then for $z\in\mathbb{D}$,
\be\label{eq-Co}
\|D_{f}(z)\|\leq\frac{4}{\pi}\frac{1}{1-|z|^{2}}.\ee

Analogy the inequality (\ref{eq-Co}), we prove the following result.

\begin{thm}\label{thm-3}
Suppose that $g\in\mathcal{C}(\overline{\mathbb{D}})$,
 $\varphi\in \mathcal{C}(\mathbb{T})$, and suppose that
$f\in\mathcal{C}^{4}(\mathbb{D})\cap\mathcal{C}(\overline{\mathbb{D}})$
maps $\overline{\mathbb{D}}$ into itself satisfying {\rm
(\ref{eq-ch-1.0})} and {\rm (\ref{eq-ch-1})}. Then for
$z\in\mathbb{D}$,
$$
\|D_{f}(z)\|
\leq\frac{4+\pi(1+2|z|+3|z|^{2})}{\pi(1-|z|^{2})}\|f\|_{\infty}+\Big(2|z|+\frac{4}{\pi}
\Big)\|\varphi\|_{\infty}+\frac{23}{24}\|g\|_{\infty}.
$$
\end{thm}

Let $\mathcal{A}(\ID)$  denote the set of all analytic functions $f$
in $\mathbb{D}$ satisfying the standard normalization:
$f(0)=f'(0)-1=0$. In the early 20th century, Landau \cite{L} showed
that there is a constant $r>0$, independent of elements in
$\mathcal{A}(\ID)$, such that $f(\mathbb{D})$ contains a disk of
radius $r$. Since then the Landau theorem has become an important
tool in geometric function theory with one complex variable (cf.
\cite{Br,Za}). Unfortunately, for general classes of functions,
there is no Landau type theorem (cf. \cite{GP, W}). To establish
analogs of the Landau type theorem for more general classes of
functions, it is necessary to restrict our focus on certain
subclasses (cf. \cite{AA, BS, BE, HG, HG1, CPW0, CP, CV, W}). Let us recall two Landau type Theorems for biharmonic mappings, which are as follows.

\begin{Thm}{\rm (\cite[Theorem 1]{AA})}\label{Thm-AA-1}
Suppose that $f(z)=|z|^{2}H(z)+K(z)$ is a biharmonic mapping, that is
$\Delta(\Delta f)=0$, in $\mathbb{D}$ such that $f(0)=
K(0)=J_{f}(0)-1=0$, where $H$ and $K$ are harmonic with
$\max\{|H(z)|, |K(z)|\}<M$, and $M$ is a positive
constant. Then there is a constant $\rho_{2}\in(0,1)$ such that $f$ is
univalent in $\mathbb{D}_{\rho_{2}}$, and
$f(\mathbb{D}_{\rho_{2}})$ contains a disk $\mathbb{D}_{R_{2}}$, where
$\rho_{2}$ satisfies the following equation:
$$\frac{\pi}{4M}-2M\rho_{2}-4M\frac{\rho_{2}}{(1-\rho_{2})^{2}}=0,$$
and
$R_{2}=\frac{\pi}{4M}\rho_{2}-2M(\rho_{2}^{3}+\rho_{2}^{2})/(1-\rho_{2}).$
\end{Thm}

\begin{Thm}{\rm (\cite[Theorem 2]{AA})}  \label{Thm-AA-2}
Suppose that $H$ denotes a harmonic mapping in $\mathbb{D}$ such
that $H(0)=J_{H}(0)-1=0$ and $H(z)<M$, where $M$ is a positive
constant. Then there is a constant $\rho_{3}$ such that
$F=|z|^{2}H$ is univalent in $\mathbb{D}_{\rho_{3}}$, and
$f(\mathbb{D}_{\rho_{3}})$ contains a disk $\mathbb{D}_{R_{3}}$,
where $\rho_{3}$ is the solution of the equation:
$$\frac{\pi}{4M}=4M\frac{\rho_{3}}{1-\rho_{3}}+2M\frac{\rho_{3}(2-\rho_{3})}{(1-\rho_{3})^{2}},$$ and $R_{3}=\frac{\pi}{4M}\rho_{3}^{2}-2M\rho_{3}^{4}/(1-\rho_{3})$.
\end{Thm}

For convenience, we make a notational convention: For
$g\in\mathcal{C}(\overline{\mathbb{D}})$ and $\varphi\in
\mathcal{C}(\mathbb{T})$, let
$\mathcal{BF}_{g}(\overline{\mathbb{D}})$ denote the class of all
complex-valued functions
$f\in\mathcal{C}^{4}(\mathbb{D})\cap\mathcal{C}(\overline{\mathbb{D}})$
satisfying {\rm (\ref{eq-ch-1.0})} and {\rm (\ref{eq-ch-1})} with
the normalization $f(0)=J_{f}(0)-1=0$.

As an application Theorem \ref{thm-3}, we establish the following Landau type Theorem for $f\in\mathcal{BF}_{g}(\overline{\mathbb{D}}).$
In particular, if $g\equiv0$, then $f\in\mathcal{BF}_{g}(\overline{\mathbb{D}})$ is biharmonic.
In this sense, the following result is a generalization of Theorems \Ref{Thm-AA-1} and \Ref{Thm-AA-2}.

\begin{thm}\label{thm-4} Suppose that $M_{1}>0$, $M_{2}\geq0$ and
$M_{3}\geq0$ are constants, and suppose that
$f\in\mathcal{BF}_{g}(\overline{\mathbb{D}})$
 satisfies the following conditions: $$\sup_{z\in\mathbb{D}}\{|f(z)|\}\leq M_{1},\;\; \sup_{z\in\mathbb{T}}\{|\varphi(z)|\}\leq
M_{2}\;\; \mbox{and}\;\; \sup_{z\in\mathbb{D}}\{|g(z)|\}\leq M_{3}.$$
 Then $f$ is univalent in
$\mathbb{D}_{r_{0}}$, and $f(\mathbb{D}_{r_{0}})$ contains a univalent disk
$\mathbb{D}_{R_{0}}$, where
  $r_{0}$ satisfies the
following equation:
\begin{eqnarray*}
&&\frac{1}{\frac{4}{\pi}(M_{1}+M_{2})+M_{1}+\frac{23}{24}M_{3}}-\frac{4(M_{1}+M_{2})}{\pi}\frac{r_{0}(2-r_{0})}{(1-r_{0})^{2}}-2M_{2}r_{0}\\
 &&-\frac{r_{0}^{2}}{1-r_{0}^{2}}\big(2M_{1}+\frac{4M_{2}}{\pi}\big)
 -\frac{101\|g\|_{\infty}r_{0}}{120}
 -\frac{2M_{1}r_{0}(2+2r_{0}+r_{0}^{2})}{(1-r_{0})(1-r_{0}^{2})}=0,
\end{eqnarray*}
and $R_{0}\geq M_0$ with
$$M_{0}= r_{0}\bigg(\frac{1}{\frac{8}{\pi}(M_{1}+M_{2})+2M_{1}+\frac{23}{12}M_{3}}+\frac{2M_{2}r_{0}^{2}}{3\pi(1-r_{0}^{2})}+
\frac{M_{1}r_{0}^{2}}{3(1-r_{0}^{2})}\bigg).$$
\end{thm}

\begin{rem}
\begin{enumerate}
\item Theorem \ref{thm-4} provides an explicit estimate of the radius $R_{0}$ of the univalent disk, which gives an answer to the open problem in
\cite[Remark 1.2]{CP}  under the additional assumption
$f\in\mathcal{C}^{4}(\mathbb{D})\cap\mathcal{C}(\overline{\mathbb{D}})$.
\item
In general, the Landau type theorem is not true for the mappings in
$\mathcal{BF}_{g}(\overline{\mathbb{D}})$. For example, let $g\equiv
1$ and $f_{k}(z)=kx+|z|^{4}/64+iy/k,$ where $z=x+iy\in\mathbb{D}$
and $k\in\{1,2,\ldots\}$. Then for all $k\in\{1,2,\ldots\}$, $f_{k}$
is univalent and $J_{f_{k}} (0)-1=f_{k}(0)=0$. But, there is no
absolute constant $R_{0}>0$ such that $\mathbb{D}_{R_{0}}$ belongs
to $f_{k}(\mathbb{D})$ for each $k$.
\end{enumerate}
\end{rem}

A continuous increasing function $\omega:\, [0,\infty)\rightarrow
[0,\infty)$ with $\omega(0)=0$ is called a {\it majorant} if
$\omega(t)/t$ is non-increasing for $t>0$ (cf. \cite{Dy1,Pav}).
Given a subset $\Omega$ of $\mathbb{C}$, a function
$\psi:~\Omega\rightarrow\mathbb{C}$ is said to belong to the {\it
Lipschitz space} $\Lambda_{\omega}(\Omega)$ if there is a positive
constant $L$ such that
$$\sup_{z_{1},z_{2}\in\Omega, z_{1}\neq z_{2}}\left\{\frac{|\psi(z_{1})-\psi(z_{2})|}{\omega(|z_{1}-z_{2}|)}\right\}\leq L.$$

It is well-known that the condition $\psi\in
\Lambda_{\omega}(\mathbb{T})$ is not enough to guarantee that its
harmonic extension $\mathcal{P}_{\psi}$ belongs to
$\Lambda_{\omega}(\mathbb{D})$, where $\omega(t)=t$. In fact,
$\mathcal{P}_{\psi}\in\Lambda_{\omega}(\mathbb{D})$ is Lipschitz
continuous if and only if the Hilbert transform of
$d\psi(e^{i\theta})/d\theta$ belongs to $L^{\infty}(\mathbb{T})$
\cite{Z}, where $\omega(t)=t$. In \cite{AKM}, the authors
established the following result for real harmonic mappings in the
unit ball $\IB^n$ of $\mathbb{R}^{n}$: For a boundary function which
is Lipschitz continuous, if its harmonic extension is quasiregular,
then this extension is also Lipschitz continuous. Recently, the
relationship of the Lipschitz continuity between the boundary
functions and their harmonic extensions has attracted much attention
\cite{K2,KM,LP}.

As the
last aim of this paper, we will investigate the Lipschitz continuity
of the solutions to the inhomogeneous biharmonic equation \eqref{eq-ch-1.0}. The result is as follows.

\begin{thm}\label{thm-5} Suppose that $\omega$ is a majorant and
$$\limsup_{t\rightarrow0^{+}}\frac{\omega(t)}{t}=c<\infty,$$ and
suppose that $f\in\mathcal{C}^{4}(\mathbb{D})$  satisfies the following equations:
$$\begin{cases}\displaystyle \Delta(\Delta f)=g &\mbox{ in }\, \mathbb{D},\\
\displaystyle  f_{\overline{z}}=\varphi &\mbox{ in }\, \mathbb{T},\\
\displaystyle f=0&\mbox{ in }\, \mathbb{T},
\end{cases}$$ where $g\in\mathcal{C}(\overline{\mathbb{D}})$, $\|\varphi\|_{\infty}<\infty$, $\|g\|_{\infty}<\infty$ and $\varphi_{1}(e^{it})=\varphi(e^{it})e^{-it}\in \Lambda_{\omega}(\mathbb{T})$ for $t\in[0,2\pi]$.
 Then $f\in\Lambda_{\omega}(\mathbb{D})$.
\end{thm}

The proofs of Theorems \ref{thm-1}, \ref{thm-2} and \ref{thm-3} will
be presented in Section \ref{csw-sec2}. Theorem \ref{thm-4} will be proved in Section
\ref{csw-sec3}, and
the proof of
Theorem \ref{thm-5} will be given in Section
\ref{csw-sec4}.

\section{Schwarz type lemmas for solutions to inhomogeneous biharmonic equations }\label{csw-sec2}
The main purpose of this section is to prove Theorems \ref{thm-1}, \ref{thm-2} and \ref{thm-3}. We start with a lemma which is used in the proof of Theorem \ref{thm-1}.

\begin{Lem}{\rm (\cite[Exercise 15 in Charpter 7]{Gar})}\label{lem-1}
For $z,$ $w\in\mathbb{D}$, suppose that $G(z,w)$ is the biharmonic
Green function defined as in {\rm (\ref{G-1})}. Then $G(z,w)\leq0.$
\end{Lem}

The proof of Theorem \ref{thm-1} also needs the following fact.

\begin{Thm}{\rm (cf. \cite{LP})}\label{Lem-2}
For any $z\in\mathbb{D}$, we have
$$\frac{1}{2\pi}\int_{0}^{2\pi}\frac{d\theta}{|1-ze^{i\theta}|^{2\alpha}}=\sum_{n=0}^{\infty}\left(\frac{\Gamma(n+\alpha)}{n!\Gamma(\alpha)}\right)^{2}|z|^{2n},$$
where $\alpha>0$ and $\Gamma$ denotes the Gamma function.
\end{Thm}

\subsection*{Proof of Theorem \ref{thm-1}}
By (\ref{eq-ch-3}) and (\ref{eq-pav1}), we have

\be\label{eq-ch-5}
\left|\mathcal{P}_{f}(z)-\frac{1-|z|^{2}}{1+|z|^{2}}\mathcal{P}_{f}(0)\right|\leq\frac{4}{\pi}\|\mathcal{P}_{f}\|_{\infty}\arctan|z|,
\ee

\beq\label{eq-ch-6}
\left|\frac{1-|z|^{2}}{2\pi}\int_{0}^{2\pi}f(e^{i\theta})\frac{\overline{z}e^{i\theta}}{(1-\overline{z}e^{i\theta})^{2}}d\theta\right|
&\leq&\frac{|z|}{2\pi}\int_{0}^{2\pi}|f(e^{i\theta})|P(z,e^{i\theta})d\theta\\
\nonumber &\leq&|z|\|\mathcal{P}_{|f|}\|_{\infty} \eeq and
\beq\label{eq-ch-7}
\bigg|(1-|z|^{2})\mathcal{P}_{\varphi_{1}}(z)-\frac{(1-|z|^{2})^{2}}{1+|z|^{2}}\mathcal{P}_{\varphi_{1}}(0)\bigg|
&\leq&
\frac{4(1-|z|^{2})}{\pi}\|\mathcal{P}_{\varphi_{1}}\|_{\infty}\arctan|z|.
\eeq

 Let $$\zeta=\frac{z-w}{1-\overline{z}w}.$$  Then,
$$w=\frac{z-\zeta}{1-\overline{z}\zeta},~w-z=\frac{\zeta(|z|^{2}-1)}{1-\overline{z}\zeta}$$ and $$|J_{w}(\zeta)|=\frac{(1-|z|^{2})^{2}}{|1-\overline{z}\zeta|^{4}},$$
which, together with Lemma \Ref{lem-1} and Theorem \Ref{Lem-2},
yields

\beq\label{eq-ch-8}
I_{0}&=&\left|\frac{1}{16\pi}\int_{\mathbb{D}}g(w)G(z,w)dA(w)\right|\leq\frac{\|g\|_{\infty}}{16\pi}\int_{\mathbb{D}}|G(z,w)|dA(w)\\
\nonumber
&=&\frac{\|g\|_{\infty}}{16\pi}\bigg(\int_{\mathbb{D}}(1-|z|^{2})(1-|w|^{2})dA(w)\\
\nonumber
&&-\int_{\mathbb{D}}|z-w|^{2}\log\left|\frac{1-z\overline{w}}{z-w}\right|^{2}dA(w)\bigg)\\
\nonumber &=&\frac{\|g\|_{\infty}(1-|z|^{2})}{32}-\frac{\|g\|_{\infty}(1-|z|^{2})^{4}}{8\pi}\int_{\mathbb{D}}\frac{|\zeta|^{2}\log\frac{1}{|\zeta|}}
{|1-\overline{z}\zeta|^{6}}dA(\zeta)\\
\nonumber
&=&\frac{\|g\|_{\infty}(1-|z|^{2})}{32}-\frac{\|g\|_{\infty}(1-|z|^{2})^{4}}{8\pi}\int_{0}^{1}\int_{0}^{2\pi}\frac{\rho^{3}\log\frac{1}{\rho}}{|1-\overline{z}\rho
e^{i\theta}|^{6}}d\theta d\rho\\ \nonumber&=&
\frac{\|g\|_{\infty}(1-|z|^{2})}{32}\\
\nonumber&&-\frac{\|g\|_{\infty}(1-|z|^{2})^{4}}{16}\sum_{n=0}^{\infty}(n+1)^{2}(n+2)^{2}|z|^{2n}\int_{0}^{1}\rho^{2n+3}\log\frac{1}{\rho}d\rho
 \\ \nonumber
&=&\frac{\|g\|_{\infty}(1-|z|^{2})}{32}-\frac{\|g\|_{\infty}(1-|z|^{4})}{64}\\
\nonumber&=&\frac{\|g\|_{\infty}(1-|z|^{2})^{2}}{64}. \eeq Hence, it follows from
(\ref{eq-ch-3}) and (\ref{eq-ch-5}) $\sim$ (\ref{eq-ch-8}) that

\begin{eqnarray*}
\bigg|f(z)-\frac{1-|z|^{2}}{1+|z|^{2}}\mathcal{P}_{f}(0)&+&
\frac{(1-|z|^{2})^{2}}{1+|z|^{2}}\mathcal{P}_{\varphi_{1}}(0)\bigg|\\
&\leq&\bigg|\mathcal{P}_{f}(z)-\frac{1-|z|^{2}}{1+|z|^{2}}\mathcal{P}_{f}(0)\bigg|\\
\nonumber&&+
\left|\frac{1}{2\pi}\int_{0}^{2\pi}f(e^{it})\frac{\overline{z}e^{it}(1-|z|^{2})}{(1-\overline{z}e^{it})^{2}}dt\right|\\
 &&+\left|(1-|z|^{2})\mathcal{P}_{\varphi_{1}}(z)-\frac{(1-|z|^{2})^{2}}{1+|z|^{2}}\mathcal{P}_{\varphi_{1}}(0)\right|+I_{0}\\
 &\leq&
\frac{4}{\pi}\|\mathcal{P}_{f}\|_{\infty}\arctan|z|+
|z|\|\mathcal{P}_{|f|}\|_{\infty}
\\
&&+\frac{4\|\mathcal{P}_{\varphi_{1}}\|_{\infty}}{\pi}(1-|z|^{2})\arctan|z|+\frac{\|g\|_{\infty}(1-|z|^{2})^{2}}{64},
\end{eqnarray*}
as required. \qed

\subsection*{Proof of Theorem \ref{thm-2}}

Since $\phi\in\mathcal{C}(\overline{\mathbb{D}})$ is analytic in
$\mathbb{D}$, we see that
\be\label{eq-r1}\frac{1}{2\pi}\int_{0}^{2\pi}\phi(e^{it})\frac{\overline{z}e^{it}(1-|z|^{2})}{(1-\overline{z}e^{it})^{2}}dt=
\frac{(1-|z|^{2})}{2\pi}\int_{0}^{2\pi}\phi(e^{it})\frac{\overline{z}e^{it}}{(1-\overline{z}e^{it})^{2}}dt=0.\ee
By using a generalized version of the classical Schwarz lemma (cf. \cite[Lemma 5.3]{Kra}), we have \be\label{eq-ch-1x}
|\mathcal{P}_{\phi}(z)|\leq\frac{|z|+|\mathcal{P}_{\phi}(0)|}{1+|z||\mathcal{P}_{\phi}(0)|}.\ee

Since applying  (\ref{eq-ch-3}), (\ref{eq-pav1}) and (\ref{eq-ch-8}) $\sim$ (\ref{eq-ch-1x}) leads to

\begin{eqnarray*}
|f(\eta)-f(r\eta)|&=&\bigg|f(\eta)+
\frac{(1-|z|^{2})^{2}}{1+|z|^{2}}\mathcal{P}_{\varphi_{1}}(0)-
\frac{(1-|z|^{2})^{2}}{1+|z|^{2}}\mathcal{P}_{\varphi_{1}}(0)-f(r\eta)\bigg|\\
&\geq&|f(\eta)|-\bigg|f(r\eta)+
\frac{(1-|z|^{2})^{2}}{1+|z|^{2}}\mathcal{P}_{\varphi_{1}}(0)\bigg|\\
&&-
\frac{(1-|z|^{2})^{2}}{1+|z|^{2}}\big|\mathcal{P}_{\varphi_{1}}(0)\big|\\
&\geq&1-\frac{|z|+|\mathcal{P}_{\phi}(0)|}{1+|z||\mathcal{P}_{\phi}(0)|}-
\left|\frac{1}{2\pi}\int_{0}^{2\pi}\phi(e^{it})\frac{\big(\overline{z}e^{it}(1-|z|^{2})\big)}{(1-\overline{z}e^{it})^{2}}dt\right|\\
 &&-\left|(1-|z|^{2})\mathcal{P}_{\varphi_{1}}(z)-\frac{(1-|z|^{2})^{2}}{1+|z|^{2}}\mathcal{P}_{\varphi_{1}}(0)\right|\\
 &&-\left|\frac{1}{16\pi}\int_{\mathbb{D}}g(w)G(z,w)dA(w)\right|-
\frac{(1-|z|^{2})^{2}}{1+|z|^{2}}|\mathcal{P}_{\varphi_{1}}(0)|\\
 &\geq&1-\frac{|z|+|\mathcal{P}_{\phi}(0)|}{1+|z||\mathcal{P}_{\phi}(0)|}-\frac{4\|\mathcal{P}_{\varphi_{1}}\|_{\infty}(1-|z|^{2})}{\pi}\arctan|z|\\
 &&-\frac{\|g\|_{\infty}}{64}(1-|z|^{2})^{2}-
\frac{(1-|z|^{2})^{2}}{1+|z|^{2}}|\mathcal{P}_{\varphi_{1}}(0)|,
\end{eqnarray*}
we know that
\begin{eqnarray*}
\liminf_{r\rightarrow1^{-}}\frac{|f(\eta)-f(r\eta)|}{1-r}&\geq&\liminf_{r\rightarrow1^{-}}\frac{(1-|\mathcal{P}_{\phi}(0)|)}{1+r|\mathcal{P}_{\phi}(0)|}
-\liminf_{r\rightarrow1^{-}}\frac{4\|\mathcal{P}_{\varphi_{1}}\|_{\infty}(1+r)}{\pi}\arctan
r\\&&-\liminf_{r\rightarrow1^{-}}\frac{\|g\|_{\infty}}{64}(1-r)(1+r)^{2}\\&&-\liminf_{r\rightarrow1^{-}}
\frac{(1-r)(1+r)^{2}}{1+r^{2}}\big|\mathcal{P}_{\varphi_{1}}(0)\big|\\
&=&\frac{1-|\mathcal{P}_{\phi}(0)|}{1+|\mathcal{P}_{\phi}(0)|}-2\|\mathcal{P}_{\varphi_{1}}\|_{\infty},
\end{eqnarray*}
which shows that the inequality \eqref{Sch-1} is true. To finish the proof of this theorem, it remains to prove the sharpness part. For this, we divide the proof into two cases.

\bca Suppose $\mathcal{P}_{\phi}(0)=0$.\eca
Let $$f(z)=\beta z$$ in $\overline{\mathbb{D}}$ with $|\beta|=1$. Then we see that the inequality {\rm(\ref{Sch-1})}
is sharp, where $z\in\overline{\mathbb{D}}$ and $|\beta|=1$.

\bca Suppose $\mathcal{P}_{\phi}(0)\neq0$.\eca

For any $a\neq 0$ in $\ID$, let $\eta=-\frac{a}{|a|}$ and
$$f(z)=\frac{z-a}{1-\overline{a}z}\;\;(z\in\overline{\mathbb{D}}).$$
By elementary calculations, we obtain that
$$\lim_{r\rightarrow1^{-}}\frac{|f(\eta)-f(r\eta)|}{1-r}=\frac{1-|a|}{1+|a|}=\frac{1-|f(0)|}{1+|f(0)|}.$$
Hence, the proof of the theorem is complete. \qed
\medskip

The following result is useful for the proof of Theorem \ref{thm-3}.

\begin{Thm}{\rm (\cite[Proposition 2.4]{K2})}\label{Thm-A}
Suppose that $X$ is an open subset of $\mathbb{R}$, and $\Omega$ a measure
space. Suppose, further, that a function
$F:$ $X\times\Omega\rightarrow\mathbb{R}$ satisfies the following
conditions:
\begin{enumerate}
\item[(1)] $F(x,w)$ is a measurable function of $x$ and $w$ jointly,
and is integrable with respect to $w$ for almost every $x\in X.$
\item[(2)] For almost every $w\in\Omega$, $F(x,w)$ is an absolutely
continuous function with respect to $x$. $($This guarantees that $\partial
F(x,w)/\partial x$ exists almost everywhere.$)$
\item[(3)] $\partial F/\partial x$ is locally integrable, that is,
for all compact intervals $[a,b]$ contained in $X$:
$$\int_{a}^{b}\int_{\Omega}\left|\frac{\partial}{\partial x}F(x,w)\right|dwdx<\infty.$$
\end{enumerate}
Then$\int_{\Omega}F(x,w)dw$ is an absolutely continuous function with respect to
$x$, and for almost every $x\in X$, its derivative exists, which is
given by
$$\frac{d}{dx}\int_{\Omega}F(x,w)dw=\int_{\Omega}\frac{\partial}{\partial x}F(x,w)dw.$$
\end{Thm}

\subsection*{Proof of Theorem \ref{thm-3}} By Theorem \Ref{Thm-A}, we have

\beq\label{eq-p1}
f_{z}(z)&=&[\mathcal{P}_{f}(z)]_{z}-\frac{1}{2\pi}\int_{0}^{2\pi}\frac{f(e^{it})\overline{z}^{2}e^{it}}{(1-\overline{z}e^{it})^{2}}dt
+\overline{z}\mathcal{P}_{\varphi_{1}}(z)\\ \nonumber
&&-(1-|z|^{2})[\mathcal{P}_{\varphi_{1}}(z)]_{z}-\frac{1}{16\pi}\int_{\mathbb{D}}g(w)G_{z}(z,w)dA(w)
\eeq and
\beq\label{eq-p2}
f_{\overline{z}}(z)&=&[\mathcal{P}_{f}(z)]_{\overline{z}}-\frac{1}{2\pi}\int_{0}^{2\pi}\frac{f(e^{it})|z|^{2}e^{it}}{(1-\overline{z}e^{it})^{2}}dt
+z\mathcal{P}_{\varphi_{1}}(z)\\ \nonumber
&&-(1-|z|^{2})[\mathcal{P}_{\varphi_{1}}(z)]_{\overline{z}}-\frac{1}{16\pi}\int_{\mathbb{D}}g(w)G_{\overline{z}}(z,w)dA(w)\\
\nonumber
&&+\frac{(1-|z|^{2})}{2\pi}\int_{0}^{2\pi}f(e^{it})\frac{e^{it}(1+\overline{z}e^{it})}{(1-\overline{z}e^{it})^{3}}dt,\eeq
where $\varphi_{1}$ is defined in {\rm (\ref{eq-ch-3})}.

By using \cite[Lemma 2.5]{LP}, we get
\be\label{eq-ch9}
\frac{\|g\|_{\infty}}{8\pi}\int_{\mathbb{D}}\big(|G_{z}(z,w)|+|G_{\overline{z}}(z,w)|\big)dA(w)\leq\frac{23}{24}\|g\|_{\infty}.\ee

By (\ref{eq-Co}) and (\ref{eq-p1}) $\sim$ (\ref{eq-ch9}),
we conclude that

\begin{eqnarray*}
\|D_{f}(z)\|&\leq&\|D_{\mathcal{P}_{f}}(z)\|+\frac{|z|^{2}}{\pi}\int_{0}^{2\pi}\frac{|f(e^{it})|}{|1-\overline{z}e^{it}|^{2}}dt+
(1-|z|^{2})\|D_{\mathcal{P}_{\varphi_{1}}}(z)\|\\
&+&2|z||\mathcal{P}_{\varphi_{1}}(z)|+\frac{23}{24}\|g\|_{\infty}+
\frac{(1-|z|^{2})}{2\pi}\int_{0}^{2\pi}|f(e^{it})|\frac{|1+\overline{z}e^{it}|}{|1-\overline{z}e^{it}|^{3}}dt\\
&\leq&\frac{4\|\mathcal{P}_{f}\|_{\infty}}{\pi(1-|z|^{2})}+\frac{2|z|^{2}\|f\|_{\infty}}{1-|z|^{2}}+\frac{4}{\pi}\|\mathcal{P}_{\varphi_{1}}\|_{\infty}
+2|z||\mathcal{P}_{\varphi_{1}}(z)|+\frac{23}{24}\|g\|_{\infty}\\
&&+\frac{1+|z|}{1-|z|}\|f\|_{\infty}\\
&\leq&\frac{4+\pi(1+2|z|+3|z|^{2})}{\pi(1-|z|^{2})}\|f\|_{\infty}+\big(
2|z|+\frac{4}{\pi}\big)\|\varphi\|_{\infty}+\frac{23}{24}\|g\|_{\infty},
\end{eqnarray*}
which is what we need. \qed

\section{A Landau type theorem for solutions to inhomogeneous biharmonic equations}\label{csw-sec3}

We will prove Theorem \ref{thm-4} in the section. First, let us recall the following result.

\begin{Thm}{\rm (\cite[Lemma 1]{CPW0})}\label{LemA}
Suppose $f$ is a harmonic mapping of $\mathbb{D}$ into $\mathbb{C}$ such
that $|f(z)|\leq M$ and
$$f(z)=\sum_{n=0}^{\infty}a_{n}z^{n}+\sum_{n=1}^{\infty}\overline{b}_{n}\overline{z}^{n}.$$
Then$|a_{0}|\leq M$ and for all $n\geq 1,$
$$|a_{n}|+|b_{n}|\leq \frac{4M}{\pi}.
$$
\end{Thm}

\subsection*{Proof of Theorem \ref{thm-4}} By elementary calculations, we have
$$G_{z}(z,w)=(\overline{z}-\overline{w})\log\left|\frac{1-z\overline{w}}{z-w}\right|^{2}+
\frac{(\overline{z}-\overline{w})(|w|^{2}-1)}{(1-z\overline{w})}+\overline{z}(1-|w|^{2})$$
and
$$G_{\overline{z}}(z,w)=(z-w)\log\left|\frac{1-z\overline{w}}{z-w}\right|^{2}+
\frac{(z-w)(|w|^{2}-1)}{(1-w\overline{z})}+z(1-|w|^{2}).$$ Then,
\beq\label{eqch-11}
G_{z}(z,w)-G_{z}(0,w)&=&2\overline{z}\log\left|\frac{1-z\overline{w}}{z-w}\right|-2\overline{w}
\big(F(z,w)-F(0,w)\big)\\
\nonumber
&&-\frac{(\overline{z}-z\overline{w}^{2})(1-|w|^{2})}{1-z\overline{w}}+\overline{z}(1-|w|^{2})
\eeq and
\begin{eqnarray*}
G_{\overline{z}}(z,w)-G_{\overline{z}}(0,w)&=&2z\log\left|\frac{1-z\overline{w}}{z-w}\right|-2w
\big(F(z,w)-F(0,w)\big)\\
&&-\frac{(z-\overline{z}w^{2})(1-|w|^{2})}{1-w\overline{z}}+z(1-|w|^{2}),
\end{eqnarray*}
 where $F(z,w)=\log\left|\frac{1-z\overline{w}}{z-w}\right|.$

By Fubini's Theorem and \cite[Inequalities (2.11) and
(2.12)]{CK}, we get

\beq\label{eqch-12} \nonumber
\int_{\mathbb{D}}\big|F(z,w)-F(0,w)\big|\frac{dA(w)}{2\pi}&\leq&\int_{\mathbb{D}}
\left(\int_{[0,z]}|F_{\varsigma}(\varsigma,w)dz+F_{\overline{\varsigma}}(\varsigma,w)d\overline{\varsigma}|\right)
\frac{dA(w)}{2\pi}\\
\nonumber
&\leq&\int_{\mathbb{D}}\left(\int_{[0,z]}\big(|F_{\varsigma}(\varsigma,w)|+|F_{\overline{\varsigma}}(\varsigma,w)|\big)|d\varsigma|\right)\frac{dA(w)}{2\pi}\\
\nonumber
&=&\int_{[0,z]}\left(\int_{\mathbb{D}}\big(|F_{\varsigma}(\varsigma,w)|+|F_{\overline{\varsigma}}(\varsigma,w)|\big)\frac{dA(w)}{2\pi}\right)|d\varsigma|\\
 &=&\int_{[0,z]}2\nu(|\varsigma|)|dz|, \eeq where
$$\nu(|\varsigma|)=\frac{1-|\varsigma|^{2}}{8|\varsigma|^{2}}
\bigg(\frac{1+|\varsigma|^{2}}{1-|\varsigma|^{2}}-\frac{1-|\varsigma|^{2}}{2|\varsigma|}\log\frac{1+|\varsigma|}{1-|\varsigma|}\bigg),$$
and $[0,z]$ denotes the segment from $0$ to $z$.

It follows from \cite[Theorem 3]{CK} that
$$\frac{1}{4}\leq\nu(|\varsigma|)\leq\frac{1}{3},$$ which, together with
(\ref{eqch-12}), gives

\be\label{eqch-13}\frac{1}{2\pi}\int_{\mathbb{D}}\big|F(z,w)-F(0,w)\big|dA(w)\leq\frac{2}{3}|z|.
\ee

By applying \cite[Inequality (2.3)]{CK}, we obtain that

\be\label{eqch-14}\int_{\mathbb{D}}\bigg|\overline{z}\log\left|\frac{1-z\overline{w}}{z-w}\right|\bigg|\frac{dA(w)}{2\pi}=\frac{|z|(1-|z|^{2})}{4}.\ee

Moreover, by elementary calculations, we have

\beq\label{eqch-15}
\int_{\mathbb{D}}\bigg|\frac{(\overline{z}-z\overline{w}^{2})(1-|w|^{2})}{1-z\overline{w}}\bigg|\frac{dA(w)}{2\pi}
&\leq&|z|\int_{\mathbb{D}}\frac{(1+|w|^{2})(1-|w|^{2})}{1-|w|}\frac{dA(w)}{2\pi}\\
\nonumber&=&\frac{77}{60}|z|. \eeq

Let
$$I_{1}=\left|\frac{1}{16\pi}\int_{\mathbb{D}}g(w)\big(G_{z}(z,w)-G_{z}(0,w)\big)dA(w)\right|$$
and
$$I_{2}=\left|\frac{1}{16\pi}\int_{\mathbb{D}}g(w)\big(G_{\overline{z}}(z,w)-G_{\overline{z}}(0,w)\big)dA(w)\right|.$$
Then by (\ref{eqch-11}) $\sim$  (\ref{eqch-15}), we deduce that
\beq\label{eqch-19}
I_{1}&\leq&\frac{\|g\|_{\infty}}{16\pi}\int_{\mathbb{D}}\big|G_{z}(z,w)-G_{z}(0,w)\big|dA(w)\\
\nonumber
&\leq&\frac{\|g\|_{\infty}|z|}{8\pi}\int_{\mathbb{D}}\bigg|\log\left|\frac{1-z\overline{w}}{z-w}\right|\bigg|dA(w)+\frac{\|g\|_{\infty}|z|}{16\pi}\int_{\mathbb{D}}(1-|w|^{2})dA(w)\\
\nonumber &&+\frac{\|g\|_{\infty}}{16\pi}\int_{\mathbb{D}}\bigg|\frac{(\overline{z}-z\overline{w}^{2})(1-|w|^{2})}{1-z\overline{w}}\bigg|dA(w)\\
\nonumber &&+
\frac{\|g\|_{\infty}}{8\pi}\int_{\mathbb{D}}\big|F(z,w)-F(0,w)\big|dA(w)\\
\nonumber
&\leq& \Big(\frac{1-|z|^{2}}{16}+\frac{43}{120}\Big)\|g\|_{\infty}|z|.
\eeq

Similarly, we can also obtain that
$$I_{2}\leq \Big(\frac{1-|z|^{2}}{16}+\frac{43}{120}\Big)\|g\|_{\infty}|z|.$$

Since $\mathcal{P}_{f}$ and $\mathcal{P}_{\varphi_{1}}$ are harmonic
in $\mathbb{D}$, we know that
$$\mathcal{P}_{f}(z)=\sum_{n=0}^{\infty}a_{n}z^{n}+\sum_{n=1}^{\infty}\overline{b}_{n}\overline{z}^{n}$$
 and
$$\mathcal{P}_{\varphi_{1}}(z)=\sum_{n=0}^{\infty}c_{n}z^{n}+\sum_{n=1}^{\infty}\overline{d}_{n}\overline{z}^{n}.$$
Thus, we infer from Theorem \Ref{LemA} that
\beq\label{eq-6r} \nonumber \big|[\mathcal{P}_{f}(z)]_{z}-
[\mathcal{P}_{f}(0)]_{z}\big|+\big|[\mathcal{P}_{f}(z)]_{\overline{z}}-
[\mathcal{P}_{f}(0)]_{\overline{z}}\big|&=&\left|\sum_{n=2}^{\infty}na_{n}z^{n-1}\right|+\left|\sum_{n=2}^{\infty}nb_{n}\overline{z}^{n-1}\right|\\
\nonumber
&\leq&\sum_{n=2}^{\infty}n\big(|a_{n}|+|b_{n}|\big)|z|^{n-1}
\\
\nonumber &\leq&\frac{4M_{1}}{\pi}\sum_{n=2}^{\infty}n|z|^{n-1}\\
&=&\frac{4M_{1}}{\pi}\frac{|z|(2-|z|)}{(1-|z|)^{2}} \eeq
and
\be\label{eq-7r} \big|[\mathcal{P}_{\varphi_{1}}(z)]_{z}-
[\mathcal{P}_{\varphi_{1}}(0)]_{z}\big|+\big|[\mathcal{P}_{\varphi_{1}}(z)]_{\overline{z}}-
[\mathcal{P}_{\varphi_{1}}(0)]_{\overline{z}}\big|\leq\frac{4M_{2}}{\pi}\frac{|z|(2-|z|)}{(1-|z|)^{2}}.\ee

Now, it follows from Theorem \Ref{Thm-A}, (\ref{eq-p1}), (\ref{eq-p2}) and
(\ref{eqch-19}) that
\beq\label{eq-8r}
|f_{z}(z)-f_{z}(0)|&=&\bigg|[\mathcal{P}_{f}(z)]_{z}-[\mathcal{P}_{f}(0)]_{z}-\frac{1}{2\pi}\int_{0}^{2\pi}\frac{f(e^{it})\overline{z}^{2}e^{it}}{(1-\overline{z}e^{it})^{2}}dt
\\ \nonumber
&&+\overline{z}\mathcal{P}_{\varphi_{1}}(z)+|z|^{2}[\mathcal{P}_{\varphi_{1}}(z)]_{z}-\big\{[\mathcal{P}_{\varphi_{1}}(z)]_{z}-[\mathcal{P}_{\varphi_{1}}(0)]_{z}\big\}\\
\nonumber
&&-\frac{1}{16\pi}\int_{\mathbb{D}}g(w)\big(G_{z}(z,w)-G_{z}(0,w)\big)dA(w)\bigg|\\
\nonumber
&\leq&\big|[\mathcal{P}_{f}(z)]_{z}-[\mathcal{P}_{f}(0)]_{z}\big|+\frac{|z|^{2}M_{1}}{1-|z|^{2}}+M_{2}|z|\\
\nonumber
&&+|z|^{2}\big|[\mathcal{P}_{\varphi_{1}}(z)]_{z}\big|+\big|[\mathcal{P}_{\varphi_{1}}(z)]_{z}-[\mathcal{P}_{\varphi_{1}}(0)]_{z}\big|+I_{1}
\eeq and

\beq\label{eq-9r}
|f_{\overline{z}}(z)-f_{\overline{z}}(0)|&=&\bigg|[\mathcal{P}_{f}(z)]_{\overline{z}}-[\mathcal{P}_{f}(0)]_{\overline{z}}-
\frac{1}{2\pi}\int_{0}^{2\pi}\frac{f(e^{it})|z|^{2}e^{it}}{(1-\overline{z}e^{it})^{2}}dt
\\ \nonumber
&&+z\mathcal{P}_{\varphi_{1}}(z)+|z|^{2}[\mathcal{P}_{\varphi_{1}}(z)]_{\overline{z}}-\big\{[\mathcal{P}_{\varphi_{1}}(z)]_{\overline{z}}-
[\mathcal{P}_{\varphi_{1}}(0)]_{\overline{z}}\big\}\\
\nonumber
&&-\frac{1}{16\pi}\int_{\mathbb{D}}g(w)\big(G_{\overline{z}}(z,w)-G_{\overline{z}}(0,w)\big)dA(w)\\
\nonumber&&+\frac{1}{2\pi}\int_{0}^{2\pi}\frac{e^{it}f(e^{it})\big(4\overline{z}e^{it}-3\overline{z}^{2}e^{2it}+\overline{z}^{3}e^{3it}\big)}{(1-\overline{z}e^{it})^{3}}dt\\
\nonumber
&&-\frac{|z|^{2}}{2\pi}\int_{0}^{2\pi}\frac{e^{it}f(e^{it})(1+\overline{z}e^{it})}{(1-\overline{z}e^{it})^{3}}dt\bigg|\\
\nonumber
&\leq&\big|[\mathcal{P}_{f}(z)]_{\overline{z}}-[\mathcal{P}_{f}(0)]_{\overline{z}}\big|+\frac{|z|^{2}M_{1}}{1-|z|^{2}}+M_{2}|z|\\
\nonumber&&+|z|^{2}\big|[\mathcal{P}_{\varphi_{1}}(z)]_{\overline{z}}\big|+\big|[\mathcal{P}_{\varphi_{1}}(z)]_{\overline{z}}-
[\mathcal{P}_{\varphi_{1}}(0)]_{\overline{z}}\big|+I_{2}\\
\nonumber&&+\frac{M_{1}|z|(4+3|z|+|z|^{2})}{(1-|z|)(1-|z|^{2})}+\frac{M_{1}|z|^{2}(1+|z|)}{(1-|z|)(1-|z|^{2})},
\eeq which, together with (\ref{eq-Co}), (\ref{eq-6r}),
(\ref{eq-7r}) and the fact $\frac{1}{2\pi}\int_{0}^{2\pi}P(z,e^{it})dt=1$,
yields

\beq\label{eq-10r}
|f_{z}(z)-f_{z}(0)|+|f_{\overline{z}}(z)-f_{\overline{z}}(0)|&\leq&\frac{4(M_{1}+M_{2})}{\pi}\frac{|z|(2-|z|)}{(1-|z|)^{2}}\\
\nonumber
&&+\frac{2|z|^{2}M_{1}}{1-|z|^{2}}+2M_{2}|z|+\frac{4M_{2}|z|^{2}}{\pi(1-|z|^{2})}\\
\nonumber &&+\frac{\|g\|_{\infty}|z|(1-|z|^{2})}{8}+\frac{43\|g\|_{\infty}|z|}{60}\\
\nonumber&&+\frac{2M_{1}|z|(2+2|z|+|z|^{2})}{(1-|z|)(1-|z|^{2})}\\
\nonumber &\leq&\tau(|z|), \eeq where the function

\begin{eqnarray*}
\tau(|z|)&=&\frac{4(M_{1}+M_{2})}{\pi}\frac{|z|(2-|z|)}{(1-|z|)^{2}}+\frac{2|z|^{2}M_{1}}{1-|z|^{2}}+2M_{2}|z|+\frac{4M_{2}|z|^{2}}{\pi(1-|z|^{2})}\\
\nonumber
&&+\frac{\|g\|_{\infty}|z|}{8}+\frac{43\|g\|_{\infty}|z|}{60}+\frac{2M_{1}|z|(2+2|z|+|z|^{2})}{(1-|z|)(1-|z|^{2})}
\end{eqnarray*}  is strictly increasing with respect to $|z|$ in $[0,1)$.

 By Theorem \ref{thm-3}, we obtain that
$$
1=J_{f}(0)=\|D_{f}(0)\|\lambda(D_{f}(0))\leq
\lambda(D_{f}(0))\Big(\frac{4}{\pi}(M_{1}+M_{2})+M_{1}+\frac{23}{24}M_{3}\Big),$$
which gives
\be\label{eq-11r}
\lambda(D_{f}(0))\geq\frac{1}{\frac{4}{\pi}(M_{1}+M_{2})+M_{1}+\frac{23}{24}M_{3}}.\ee

To prove the univalence of $f$ in $\mathbb{D}_{r_{0}}$, we
choose two points  $z_{1}\not=z_{2}\in\mathbb{D}_{r_{0}}$, where $r_{0}$ satisfies the
following equation:
\beq\label{eq-1e}
&&\frac{1}{\frac{4}{\pi}(M_{1}+M_{2})+M_{1}+\frac{23}{24}M_{3}}-\frac{4(M_{1}+M_{2})}{\pi}\frac{r_{0}(2-r_{0})}{(1-r_{0})^{2}}-2M_{2}r_{0}\\
\nonumber
 &&-\frac{r_{0}^{2}}{1-r_{0}^{2}}\big(2M_{1}+\frac{4M_{2}}{\pi}\big)
 -\frac{101\|g\|_{\infty}r_{0}}{120}
 -\frac{2M_{1}r_{0}(2+2r_{0}+r_{0}^{2})}{(1-r_{0})(1-r_{0}^{2})}=0.
 \eeq
Then(\ref{eq-10r}) and (\ref{eq-11r}) guarantee that
\begin{eqnarray*}
|f(z_{2})-f(z_{1})|&=&\left|\int_{[z_{1},z_{2}]}f_{z}(z)dz+f_{\overline{z}}(z)d\overline{z}\right|\\
&\geq&\left|\int_{[z_{1},z_{2}]}f_{z}(0)dz+f_{\overline{z}}(0)d\overline{z}\right|\\
&&-\left|\int_{[z_{1},z_{2}]}\big(f_{z}(z)-f_{z}(0)\big)dz+\big(f_{\overline{z}}(z)-f_{\overline{z}}(0)\big)d\overline{z}\right|\\
&\geq&\lambda(D_{f}(0))|z_{2}-z_{1}|\\
&&-\int_{[z_{1},z_{2}]}\big(|f_{z}(z)-f_{z}(0)|+|f_{\overline{z}}(z)-f_{\overline{z}}(0)|\big)|dz|\\
&>&|z_{2}-z_{1}|\bigg(\frac{1}{\frac{4}{\pi}(M_{1}+M_{2})+M_{1}+\frac{23}{24}M_{3}}-\frac{4M_{2}r_{0}^{2}}{\pi(1-r_{0}^{2})}\\
 &&-\frac{4(M_{1}+M_{2})}{\pi}\frac{r_{0}(2-r_{0})}{(1-r_{0})^{2}}-\frac{2r_{0}^{2}M_{1}}{1-r_{0}^{2}}-2M_{2}r_{0}\\
 &&-\frac{\|g\|_{\infty}r_{0}}{8}-\frac{43\|g\|_{\infty}r_{0}}{60}-\frac{2M_{1}r_{0}(2+2r_{0}+r_{0}^{2})}{(1-r_{0})(1-r_{0}^{2})}\bigg)\\
 &=&0,
\end{eqnarray*} which implies that $f(z_{2})\neq f(z_{1})$. Thus, from the arbitrariness of $z_{1}$ and $z_{2}$, the univalence of $f$ follows.

To finish the proof of this theorem, it remains to show that the image $f(\mathbb{D}_{r_{0}})$ contains a disk. To reach this goal, let $\zeta=r_{0}e^{i\theta}\in\partial\mathbb{D}_{r_{0}}$. Then we infer from (\ref{eq-10r}) and (\ref{eq-11r}) that
\begin{eqnarray*}
|f(\zeta)-f(0)|&=&\left|\int_{[0,\zeta]}f_{z}(z)dz+f_{\overline{z}}(z)d\overline{z}\right|\\
&\geq&\left|\int_{[0,\zeta]}f_{z}(0)dz+f_{\overline{z}}(0)d\overline{z}\right|\\
&&-\left|\int_{[0,\zeta]}\big(f_{z}(z)-f_{z}(0)\big)dz+\big(f_{\overline{z}}(z)-f_{\overline{z}}(0)\big)d\overline{z}\right|\\
&\geq&\lambda(D_{f}(0))r_{0}-\int_{[0,\zeta]}\big(|f_{z}(z)-f_{z}(0)|+|f_{\overline{z}}(z)-f_{\overline{z}}(0)|\big)|dz|\\
&\geq&\lambda(D_{f}(0))r_{0}-\bigg(\frac{4(M_{1}+M_{2})}{\pi}\frac{(2-r_{0})}{(1-r_{0})^{2}}\int_{[0,\zeta]}|z||dz|\\
&&+\frac{2M_{1}}{1-r_{0}^{2}}\int_{[0,\zeta]}|z|^{2}|dz|+2M_{2}\int_{[0,\zeta]}|z||dz|\\
\nonumber
&&+\frac{4M_{2}}{\pi(1-r_{0}^{2})}\int_{[0,\zeta]}|z|^{2}|dz|+\frac{101\|g\|_{\infty}}{120}\int_{[0,\zeta]}|z||dz|\\
 &&+\frac{2M_{1}(2+2r_{0}+r_{0}^{2})}{(1-r_{0})(1-r_{0}^{2})}\int_{[0,\zeta]}|z||dz|\bigg)\\
&=&r_{0}\bigg(\frac{1}{\frac{4}{\pi}(M_{1}+M_{2})+M_{1}+\frac{23}{24}M_{3}}-\frac{2(M_{1}+M_{2})}{\pi}\frac{r_{0}(2-r_{0})}{(1-r_{0})^{2}}\\
&&-\frac{2r_{0}^{2}M_{1}}{3(1-r_{0}^{2})}-M_{2}r_{0}-\frac{4M_{2}r_{0}^{2}}{3\pi(1-r_{0}^{2})}-\frac{101\|g\|_{\infty}r_{0}}{240}\\
 &&-\frac{M_{1}r_{0}(2+2r_{0}+r_{0}^{2})}{(1-r_{0})(1-r_{0}^{2})}\bigg)\\
 &=&r_{0}\bigg(\frac{1}{\frac{8}{\pi}(M_{1}+M_{2})+2M_{1}+\frac{23}{12}M_{3}}+\frac{2M_{2}r_{0}^{2}}{3\pi(1-r_{0}^{2})}+
\frac{r_{0}^{2}M_{1}}{3(1-r_{0}^{2})}\bigg),
\end{eqnarray*}
which implies that $f(\mathbb{D}_{r_{0}})$ contains a univalent disk $\mathbb{D}_{R_{0}}$ with the radius $R_0$ satisfying
$$R_{0}\geq r_{0}\bigg(\frac{1}{\frac{8}{\pi}(M_{1}+M_{2})+2M_{1}+\frac{23}{12}M_{3}}+\frac{2M_{2}r_{0}^{2}}{3\pi(1-r_{0}^{2})}+
\frac{M_{1}r_{0}^{2}}{3(1-r_{0}^{2})}\bigg).$$ Hence, the proof of this
theorem is complete. \qed

\section{Lipschitz type spaces on solutions to inhomogeneous biharmonic equations }\label{csw-sec4}

\subsection*{Proof of Theorem \ref{thm-5}} To prove Theorem \ref{thm-5}, let $z=re^{i\theta}\in\mathbb{D}$, where $r\in[0,1)$ and
$\theta\in[0,2\pi]$. Then by (\ref{eq-ch-3}), we have

\beq\label{eqL1}
f(z)&=&-(1-|z|^{2})\mathcal{P}_{\varphi_{1}}(z)-\frac{1}{16\pi}\int_{\mathbb{D}}g(w)G(z,w)dA(w)\\
\nonumber
&=&-\frac{(1-|z|^{2})}{2\pi}\int_{0}^{2\pi}P(z,e^{it})\big(\varphi_{1}(e^{it})-\varphi_{1}(e^{i\theta})\big)dt\\
\nonumber
&&-(1-|z|^{2})\varphi_{1}(e^{i\theta})-\frac{1}{16\pi}\int_{\mathbb{D}}g(w)G(z,w)dA(w),
\eeq which, together with Theorem \Ref{Thm-A}, gives

\begin{eqnarray*}
f_{z}(z)&=&-\frac{(1-|z|^{2})}{2\pi}\int_{0}^{2\pi}P_{z}(z,e^{it})\big(\varphi_{1}(e^{it})-\varphi_{1}(e^{i\theta})\big)dt\\
&&+\frac{\overline{z}}{2\pi}\int_{0}^{2\pi}P(z,e^{it})\big(\varphi_{1}(e^{it})-\varphi_{1}(e^{i\theta})\big)dt\\
&&+\overline{z}\varphi_{1}(e^{i\theta})-\frac{1}{16\pi}\int_{\mathbb{D}}g(w)G_{z}(z,w)dA(w)
\end{eqnarray*} and

\begin{eqnarray*}
f_{\overline{z}}(z)&=&-\frac{(1-|z|^{2})}{2\pi}\int_{0}^{2\pi}P_{\overline{z}}(z,e^{it})\big(\varphi_{1}(e^{it})-\varphi_{1}(e^{i\theta})\big)dt\\
&&+\frac{z}{2\pi}\int_{0}^{2\pi}P(z,e^{it})\big(\varphi_{1}(e^{it})-\varphi_{1}(e^{i\theta})\big)dt\\
&&+z\varphi_{1}(e^{i\theta})-\frac{1}{16\pi}\int_{\mathbb{D}}g(w)G_{\overline{z}}(z,w)dA(w).
\end{eqnarray*}

Let
$$E_{1}(t)=\{t\in[0,2\pi]:~|e^{it}-e^{i\theta}|\leq1-r\}\;\; \mbox{and}\;\; E_{2}(t)=\{t\in[0,2\pi]:~|e^{it}-e^{i\theta}|>1-r\}.$$ Then,
\be\label{eq-L}|e^{it}-e^{i\theta}|\leq|e^{it}-z|+|z-e^{i\theta}|=|e^{it}-z|+1-r\leq2|e^{it}-z|.\ee

Since $\varphi_{1}\in\Lambda_{\omega}(\mathbb{T})$, there is  a
positive constant $L$ such that
\be\label{eq-L2}|\varphi_{1}(e^{it})-\varphi_{1}(e^{i\theta})|\leq
L\omega(|e^{it}-e^{i\theta}|),\ee
and so, we deduce from (\ref{eq-L}) and (\ref{eq-L2}) that
\beq\label{eqL3}
|I_{3}(z)|&=&\left|\frac{1}{2\pi}\int_{0}^{2\pi}P_{z}(z,e^{it})\big(\varphi_{1}(e^{it})-\varphi_{1}(e^{i\theta})\big)dt\right|\\
\nonumber
&\leq&\frac{1}{2\pi}\int_{0}^{2\pi}\frac{1}{|1-ze^{-it}|^{2}}\big|\varphi_{1}(e^{it})-\varphi_{1}(e^{i\theta})\big|dt\\
\nonumber
&=&\frac{1}{2\pi}\int_{E_{1}(t)}\frac{1}{|1-ze^{-it}|^{2}}\big|\varphi_{1}(e^{it})-\varphi_{1}(e^{i\theta})\big|dt\\
\nonumber
&&+\frac{1}{2\pi}\int_{E_{2}(t)}\frac{1}{|1-ze^{-it}|^{2}}\big|\varphi_{1}(e^{it})-\varphi_{1}(e^{i\theta})\big|dt\\
\nonumber
&\leq&\frac{L}{2\pi}\int_{E_{1}(t)}\frac{|e^{it}-e^{i\theta}|}{|1-ze^{-it}|^{2}}\frac{\omega(|e^{it}-e^{i\theta}|)}{|e^{it}-e^{i\theta}|}dt\\
\nonumber
&&+\frac{L}{2\pi}\int_{E_{2}(t)}\frac{|e^{it}-e^{i\theta}|}{|1-ze^{-it}|^{2}}\frac{\omega(|e^{it}-e^{i\theta}|)}{|e^{it}-e^{i\theta}|}dt\\
\nonumber
&\leq&\frac{Lc}{2\pi(1-|z|)}\int_{E_{1}(t)}dt+\frac{Lc}{\pi}\int_{E_{2}(t)}\frac{|e^{it}-z|}{|1-ze^{-it}|^{2}}dt\\
\nonumber
&\leq&\frac{Lc}{(1-|z|)}\left(\frac{1}{2\pi}\int_{E_{1}(t)}dt+\frac{1}{2\pi}\int_{E_{2}(t)}dt+\frac{1}{2\pi}\int_{E_{2}(t)}dt\right)\\
&\leq& \nonumber\frac{2Lc}{(1-|z|)} \eeq and

\beq\label{eqL4}
|I_{4}(z)|&=&\left|\frac{1}{2\pi}\int_{0}^{2\pi}P(z,e^{it})\big(\varphi_{1}(e^{it})-\varphi_{1}(e^{i\theta})\big)dt\right|\\
\nonumber
&\leq&\frac{L}{2\pi}\int_{0}^{2\pi}P(z,e^{it})\omega(|e^{it}-e^{i\theta}|)dt\\
\nonumber &\leq&2L\omega(2). \eeq
Thus, by (\ref{eqL3}) and (\ref{eqL4}), we get that
\beq\label{eqL5}|f_{z}(z)|&\leq&|I_{3}(z)|(1-|z|^{2})+|z||I_{4}(z)|+|z|\|\varphi_{1}\|_{\infty}\\
\nonumber
&&+\frac{\|g\|_{\infty}}{16\pi}\int_{\mathbb{D}}|G_{z}(z,w)|dA(w)\\
\nonumber
&\leq&4Lc+2L\omega(2)+\|\varphi_{1}\|_{\infty}+\frac{\|g\|_{\infty}}{16\pi}\int_{\mathbb{D}}|G_{z}(z,w)|dA(w).
\eeq

Similarly, we have
\be\label{eqL6}|f_{\overline{z}}(z)|\leq4Lc+2L\omega(2)+\|\varphi_{1}\|_{\infty}+\frac{\|g\|_{\infty}}{16\pi}\int_{\mathbb{D}}|G_{\overline{z}}(z,w)|dA(w).
\ee
Hence, it follows from (\ref{eq-ch9}), (\ref{eqL5}) and (\ref{eqL6}) that
\beq\label{eqL7}\nonumber
\|D_{f}(z)\|&\leq&8Lc+4L\omega(2)+2\|\varphi_{1}\|_{\infty}+\frac{\|g\|_{\infty}}{16\pi}\int_{\mathbb{D}}\big(|G_{z}(z,w)|+|G_{\overline{z}}(z,w)|\big)dA(w)\\
\nonumber
&\leq&8Lc+4L\omega(2)+2\|\varphi_{1}\|_{\infty}+\frac{23}{48}\|g\|_{\infty},
\eeq which implies
$$\sup_{z,w\in\Omega, z\neq w}\left\{\frac{|f(z)-f(w)|}{|z-w|}\right\}\leq
8Lc+4L\omega(2)+2\|\varphi_{1}\|_{\infty}+\frac{23}{48}\|g\|_{\infty},$$
as required. Thus, the proof of this theorem is complete. \qed

\bigskip

{\bf Acknowledgements:} We thank the referee for providing
constructive comments and help in improving this paper. This
research was partly supported by the National Natural Science
Foundation of China ( No. 11571216), the Science and Technology Plan
Project of Hunan Province (No. 2016TP1020), the Science and
Technology Plan Project of Hengyang City (2017KJ183), and the
Construct Program of the Key Discipline in Hunan Province.

\normalsize

\end{document}